\documentclass{amsart}
\usepackage{geometry,amssymb,color,hyperref}
\usepackage{graphicx} 

\theoremstyle{plain}

\theoremstyle{definition}

\newcommand{\IZ}{\mathbb Z}

\theoremstyle{remark}
\newtheorem{remark}{Remark}[section]

\newtheorem{conjecture}{Conjecture}[section]

\title[New $S(2,9,369)$, unitals of order $6$ and other new Steiner systems]{Existence of $S(2,9,369)$, new unitals of order $6$ and other Steiner systems with block length $\ge 7$}
\author{Ivan Hetman}
\date{November 2025}

\begin{document}

    \begin{abstract} Whereas Steiner systems $S(2,k,v)$ with block length $k \le 5$ have large amount of examples and the existence is established for all admissible $v$, for $k\ge 6$ only few examples are known even for decided cases. In this paper the existence of $S(2,9,369)$ is established and some new examples for other admissible pairs $(k,v)$ are given. In particular, lots of new unitals of order $6$ (or $S(2,7,217)$) together with $S(2,7,175)$, $S(2,7,259)$, $S(2,8,120)$, $S(2,8,504)$, $S(2,9,513)$ are presented. Found examples suggest two conjectures on infinite series of designs.
    \end{abstract}

    \maketitle

    \section{Introduction}

    In this paper different $S(2,k,v)$ with $k \ge 7$ will be provided. Common source for designs with $k \ge 7$ are difference families over various fields $\mathbb F_v$ which establish infinite series \cite{CWZ}, that then can be extended to other designs using recursive constructions \cite{Greig}. This approach however leaves undecided cases. For instance, the existence of $S(2,7,127)$ and $S(2,7,211)$ is undecided despite $127$ and $211$ are prime numbers. Therefore for generating $S(2,k,v)$ with relatively small $v$ another approaches are used.

    To generate new designs, we take a specific group $G$ action on a set $X$. The $G$-space $X$ consists of $n$ orbits $X=\bigcup_{i=1}^n X_i$. Then we need to find basic blocks $B_1,\dots,B_m$ such that $|gB_i\cap B_i|\ge 2 \implies gB_i = B_i$. Obviously the amount of possible orbit constructions explodes when we increase the number of orbits. Also when we increase number of orbits, it becomes harder or impossible to generate Steiner systems due to lack of computational power. The following construction can give us a new Steiner system. Given an existing Steiner system $X$ with large automorphism group $Aut(X)$, take a subgroup $G < Aut(X)$ of the automorphism group and find $X_1,\dots,X_n$ orbits for action of $G$ on $X$. This gives us explicit construction on which we can try to find new Steiner systems.

    For example, consider design $X$ of $S(2,7,175)$ discovered in \cite{JT}. It has large automorphism group $|Aut(X)|=4200$. Among subgroups of $Aut(X)$ there is a group $G\cong (\IZ_5 \times \IZ_5) \rtimes \IZ_6$ which acts on points of this design producing three orbits of size $75,75,25$. Beginning from this construction, it is possible to generate $4$ new designs together with $2$ known.

    Additionally, one can consider constructions for smaller $v$ and $k$ and try to transform those constructions to bigger $v$ and $k$ by using groups similar to generating group for smaller values. Using this technique, new $S(2,7,175)$ were found and the undecided case $S(2,9,369)$ was eventually closed. This also allows to pose two conjectures on the existence of two infinite series of Steiner systems.

    \section{New Steiner systems}

    In this section new Steiner systems are presented. Generally there is no convenient way to describe Steiner system using convenient notation, therefore next way of presenting is chosen. All blocks of corresponding systems can be found at URL \url{https://github.com/Ihromant/math-utils/blob/master/src/test/resources/ref/}. Each line consists of list of $0$-based blocks of the designs preceded by fingerprint. The fingerprint was introduced in \cite{Het} and can be used to distinguish non-isomorphic designs. If two designs have different finerprints, then they are non-isomorphic. Table below shows information about groups that were used to generate designs together with respecting orbit structure. My source of small groups is GAP \cite{GAP4}. The following notation will be used to describe the orbit structure. For instance, $168+1\times 7$ means that we have one large orbit of size $168$ and seven orbits of size $1$. GAP ID {\tt (168,42)} means that we are considering {\tt SmallGroup(168,42)}. $\ge x$ in \# of designs means that search was unexhaustive for this construction.
    \begin{center}
        \begin{tabular}{|c c c c c c c|}
            \hline
            {$v$} & {$k$} & {GAP ID} & {group structure} & {orbits} & {\# of des} & {Known} \\
            \hline
            $175$ & $7$ & {\tt (150,6)} & {\tt (C5\,x\,C5)\,:\,C6} & $75\times 2 + 25$ & $6$ & 2 \cite{JT} \cite{Krc} \\
            $175$ & $7$ & {\tt (168,42)} & {\tt PSL(3,2)} & $168 + 1\times 7$ & $1$ & - \\
            $175$ & $7$ & {\tt (168,43)} & {\tt (C2\,x\,C2\,x\,C2)\,:\,(C7\,:\,C3)} & $168 + 1\times 7$ & $2$ & - \\
            $175$ & $7$ & {\tt (168,51)} & {\tt C2\,x\,C2\,x\,C2\,x\,(C7\,:\,C3)} & $168 + 1\times 7$ & $1$ & - \\
            $175$ & $7$ & {\tt (504,74)} & {\tt (C7\,:\,C3)\,x\,SL(2,3)} & $168 + 7$ & $1$ & - \\
            $175$ & $7$ & {\tt (504,127)} & {\tt C3\,x\,(Q8\,:\,(C7\,:\,C3))} & $168 + 7$ & $3$ & - \\
            $175$ & $7$ & {\tt (504,157)} & {\tt C3\,x\,PSL(3,2)} & $168 + 7$ & $2$ & - \\
            $175$ & $7$ & {\tt (504,158)} & {\tt C3\,x\,((C2\,x\,C2\,x\,C2)\,:\,(C7\,:\,C3))} & $168 + 7$ & $4$ & - \\
            $175$ & $7$ & {\tt (504,159)} & {\tt S4\,x\,(C7\,:\,C3)} & $168 + 7$ & $1$ & - \\
            $175$ & $7$ & {\tt (504,163)} & {\tt C2\,x\,((C7\,:\,C3)\,x\,A4)} & $168 + 7$ & $1$ & - \\
            $217$ & $7$ & {\tt (648,93)} & {\tt Q8\,:\,((C9\,x\,C3)\,:\,C3)} & $216 + 1$ & $3$ & - \\
            $217$ & $7$ & {\tt (648,94)} & {\tt Q8\,:\,((C9\,x\,C3)\,:\,C3)} & $216 + 1$ & $1$ & - \\
            $217$ & $7$ & {\tt (648,98)} & {\tt Q8\,:\,((C9\,x\,C3)\,:\,C3)} & $216 + 1$ & $1$ & - \\
            $217$ & $7$ & {\tt (648,101)} & {\tt Q8\,:\,((C9\,x\,C3)\,:\,C3)} & $216 + 1$ & $9$ & - \\
            $217$ & $7$ & {\tt (648,102)} & {\tt Q8\,:\,((C3\,x\,C3\,x\,C3)\,:\,C3)} & $216 + 1$ & $5$ & - \\
            $217$ & $7$ & {\tt (648,104)} & {\tt Q8\,:\,((C3\,x\,C3\,x\,C3)\,:\,C3)} & $216 + 1$ & $1$ & - \\
            $259$ & $7$ & {\tt (777,3)} & {\tt C7\,:\,(C37\,:\,C3)} & $259$ & $46$ & 1 \cite{Abel} \\
            $259$ & $7$ & {\tt (777,4)} & {\tt C7\,:\,(C37\,:\,C3)} & $259$ & $46$ & - \\
            \hline
            $120$ & $8$ & {\tt (120,34)} & {\tt S5} & $120$ & $3$ & - \\
            $120$ & $8$ & {\tt (112,41)} & {\tt C2\,x\,((C2\,x\,C2\,x\,C2)\,:\,C7)} & $112 + 7 + 1$ & $3$ & - \\
            $120$ & $8$ & {\tt (112,41)} & {\tt C2\,x\,((C2\,x\,C2\,x\,C2)\,:\,C7)} & $112 + 8$ & $3$ & - \\
            $512$ & $8$ & - & {\tt SL(2,8)\,x\,C7} & $504 + 7 + 1$ & $\ge 18$ & - \\
            \hline
            $369$ & $9$ & {\tt (720,766)} & {\tt C2\,x\,A6} & $360 + 2 \times 4 + 1$ & $\ge 10$ & - \\
            $513$ & $9$ & - & {\tt SL(2,8)\,x\,C9} & $504 + 9$ & $\ge 6$ & - \\
            \hline
        \end{tabular}
    \end{center}

    \begin{remark} The description of the $S(2,7,259)$ in \cite{HoCD} is not reproducible. Correct description is present in original paper of Abel \cite{Abel}.
    \end{remark}
    \begin{remark} Known $4$ examples of $S(2,7,217)$ \cite{Krc} are generated by cyclic difference families while new designs are $1$-rotational, so they are different.
    \end{remark}
    \begin{remark} Known examples of $S(2,8,120)$ embeddable into projective planes of order $16$ are listed in \cite{Pent}, but their automorphism group orders don't match new designs automorphism group orders. Therefore new $S(2,8,120)$ are not embeddable in known projective planes of order $16$.
    \end{remark}
    \begin{remark} The new $S(2,9,369)$ resolves one of undecided cases listed in \cite[II.3.12]{HoCD} with $t=5$. By applying recursive construction, from any $S(2,9,369)$ design, we can obtain $S(2,9,369\cdot 9)=S(2,9,3321)=S(2,9,46\cdot 72+9)$ design, which also belongs to the list of undecided cases with $t=46$.
    \end{remark}

    \section{Infinite series conjecture}

    We suggest that the found designs begin two infinite design series.
    \begin{conjecture} For any finite field $\mathbb F_k$ and the group $G=PSL(2,\mathbb F_k)$ there exists a Steiner system $S(2,k,|G| + k)$ based on a $G$-space that splits into one orbit of size $|G|$ and $k$ orbits of size $1$. The conjecture was confirmed for fields of cardinality $k \le 9$.
    \end{conjecture}

    \begin{conjecture} For any finite field $\mathbb F_k$ and the group $G=SL(2,\mathbb F_k)$ there exists a Steiner system (and unital) $S(2,k + 1,k^3 + 1)$ based on a $G$-space that splits into one orbit of size $|G|$ and $k + 1$ orbits of size $1$. The conjecture was confirmed for fields of cardinality $k\in\{2, 3, 4, 5, 8\}$. These designs can correspond to the unitals in Hall planes described in \cite{Gr}. For $k\in \{3, 4\}$ the corresponding unitals are isomorphic to Hall plane unitals.
    \end{conjecture}

    \section{Acknowledgement}

    I would like to express my sincere thanks to Taras Banakh and Alex Ravsky for various suggestions that allowed me to adapt my algorithm for large groups and orbit sizes which eventually led to newly discovered Steiner systems.

\end{document}